# A CHAOTIC REPRESENTATION PROPERTY OF THE MULTIDIMENSIONAL DUNKL PROCESSES

By Léonard Gallardo and Marc Yor

*Université de Tours and Université Paris 6 et IUF*

Dunkl processes are martingales as well as càdlàg homogeneous Markov processes taking values in $\mathbb{R}^d$ and they are naturally associated with a root system. In this paper we study the jumps of these processes, we describe precisely their martingale decompositions into continuous and purely discontinuous parts and we obtain a Wiener chaos decomposition of the corresponding $L^2$ spaces of these processes in terms of adequate mixed multiple stochastic integrals.

**1. Introduction.**

1.1. *Defining the Dunkl processes.* Various motivations, originating either from physical problems, or purely mathematical ones may lead to the study of the class of $\mathbb{R}^d$-valued Feller Markov processes $\{(X_t)_{t \geq 0}; (\mathbb{P}_x)_{x \in \mathbb{R}^d}\}$ which, moreover, satisfy the properties:

  (i) under $\mathbb{P}_x(X_t, t \geq 0)$ is a martingale;
  (ii) the process $(X_t, t \geq 0)$ satisfies the Brownian scaling property, that is,

$$\forall c > 0 \qquad \{(X_{ct}, t \geq 0); \mathbb{P}_x\} \stackrel{(d)}{=} \{(\sqrt{c} X_t, t \geq 0); \mathbb{P}_{x/\sqrt{c}}\}.$$

Lamperti [15] gives a striking description of all $\mathbb{R}_+$-valued Feller–Markov processes with the scaling property in terms of Lévy processes. This may motivate studies of multidimensional Markov processes with the scaling property; in this paper, we are interested in the particular class of such processes, obtained with the so-called *Dunkl processes*, which have already been studied by Rösler [18], Rösler and Voit [20], and more recently, by Gallardo and Yor [11, 12].

The Dunkl Markov processes in $\mathbb{R}^d$ are càdlàg Markov processes with infinitesimal generators:

$$\mathcal{L}_k = \tfrac{1}{2} L_k \equiv \tfrac{1}{2} \sum_{i=1}^d T_i^2,$$

where, for $1 \leq i \leq d$, $T_i$ is the differential–difference operator defined for $u \in C^1(\mathbb{R}^d)$, by

(1) $$T_i u(x) = \frac{\partial u(x)}{\partial x_i} + \sum_{\alpha \in R_+} k(\alpha) \alpha_i \frac{u(x) - u(\sigma_\alpha x)}{\langle \alpha, x \rangle},$$

where $\langle \cdot, \cdot \rangle$ is the Euclidean scalar product, $R$ is a root system in $\mathbb{R}^d$, $R_+$ is a positive subsystem, $k$ is a nonnegative multiplicity function defined on $R$ and invariant by the finite reflection group $W$ associated with $R$ and $\sigma_\alpha$ is the reflection with respect to the hyperplane $H_\alpha$ orthogonal to $\alpha$ [7, 8, 9]. For later convenience, we will suppose that the roots are normalized as follows:

(2) $$\forall \alpha \in R \qquad \langle \alpha, \alpha \rangle = |\alpha|^2 = 2.$$

As a consequence, the reflection $\sigma_\alpha$ is given by

(3) $$\sigma_\alpha(x) = x - \langle \alpha, x \rangle \alpha, \qquad x \in \mathbb{R}^d.$$

Computed for $u \in C_c^2(\mathbb{R}^d)$, the infinitesimal generator $\mathcal{L}_k$ of the Dunkl process $X = (X_t)_{t \geq 0}$ with multiplicity function $k$ reads

(4) $$\begin{aligned}\mathcal{L}_k u(x) &= \frac{1}{2} L_k u(x) \\ &= \frac{1}{2} \triangle u(x) + \sum_{\alpha \in R_+} k(\alpha) \left[ \frac{\langle \nabla u(x), \alpha \rangle}{\langle \alpha, x \rangle} + \frac{u(\sigma_\alpha x) - u(x)}{\langle \alpha, x \rangle^2} \right],\end{aligned}$$

where $\Delta$ is the Laplacian operator on $\mathbb{R}^d$.

1.2. *On the Dunkl theory.* In Physics, the Dunkl operators $T_i$ and other differential–reflection operators, appear as a crucial tool to study the integrability of the Hamiltonian operators of some Calogero–Moser–Sutherland quantum mechanical systems. Precisely those systems where the particles are distributed on a line or a torus and interact pairwisely through long range potentials of inverse square type [6]. For example, the Dunkl operator $L_k$ appears as a gauge transform of the Calogero Hamiltonian with harmonic confinement (see [19] where the reader will find a nice introduction to the analytical aspects of these questions). In Analysis, the theory initiated by Dunkl has allowed to develop a detailed spectral study of these operators (see [9] and many references therein).

1.3. *Main results in this paper.* The aim of this paper is to show that a Wiener chaos decomposition property holds for the multidimensional Dunkl process $X$ which is also a martingale. Let us give a short description of our results.

The starting point of our study is the martingale decomposition of $X$ into its continuous part which is a $d$-dimensional Brownian motion $\{B_t = (B_t^{(1)}, \ldots, B_t^{(d)}), t \geq 0\}$ and its purely discontinuous part $(\eta_t)_{t \geq 0}$ which can be explicitly written as a sum,

$$
(5) \qquad \eta_t = \sum_{\alpha \in R_+} \sqrt{k(\alpha)} M_t^\alpha \alpha,
$$

where $(M_t^\alpha, t \geq 0)$ ($\alpha \in R_+$) are one-dimensional orthogonal normal martingales (i.e., $\langle M^\alpha, M^\alpha \rangle_t = t$ and $[M^\alpha, M^\beta]_t = 0$ for $\alpha \neq \beta$). The martingale $(M_t^\alpha)$ is precisely shown to be the compensated sum of the jumps performed by $(X_t)$ in the direction of the root $\alpha$.

Now, we let $N = d + \mathrm{card}(R_+)$, we define a numbering of the roots $\alpha \in R_+$ of the form $\alpha_{d+1}, \ldots, \alpha_N$ and for any integer $\varepsilon \in [1, N]$, we consider the real-valued martingale $Z^\varepsilon = (Z_t^\varepsilon)_{t \geq 0}$ given by

$$
Z_t^\varepsilon = \begin{cases} B_t^{(\varepsilon)}, & \text{if } \varepsilon \in [1, d], \\ M_t^{\alpha_\varepsilon}, & \text{if } \varepsilon \in [d+1, N]. \end{cases}
$$

The Wiener chaos $C_n(X)$ of order $n \geq 1$ is then defined as the subspace of $L^2(\mathcal{F}_\infty)$—the $L^2$ space of the $\sigma$ field generated by the Dunkl process $X$—spanned by the "mixed" stochastic integrals

$$
\int_0^\infty dZ_{u_1}^{\varepsilon_1} \int_0^{(u_1)-} dZ_{u_2}^{\varepsilon_2} \cdots \int_0^{(u_{n-1})-} dZ_{u_n}^{\varepsilon_n} f_n(u_1, \ldots, u_n),
$$

for $(\varepsilon_1, \ldots, \varepsilon_n) \in [1, N]^n$ and $f_n : \Delta_n \to \mathbb{R}$ a measurable function such that

$$
\int_{\Delta_n} du_1 \cdots du_n\, f_n^2(u_1, \ldots, u_n) < \infty,
$$

where

$$
\Delta_n = \{(u_1, \ldots, u_n); 0 < u_n < u_{n-1} < \cdots < u_1\}.
$$

Finally we can state the main result of our paper in the form of the direct sum decomposition:

$$
L^2(\mathcal{F}_\infty) = \bigoplus_{n=0}^\infty C_n(X),
$$

where $C_0(X)$ is the subspace of constants.

For the reader's convenience, we emphasize that an extensive investigation of the one-dimensional case has been performed in [11] which may be considered as a guide for the present more general study. Precisely, our two main results in [11] are the following:

(a) A skew-product decomposition of the one-dimensional Dunkl process $(X_t)_{t\geq 0}$ as

$$X_t = |X_t|(-1)^{\mathcal{N}_{A_t}}, \qquad t \geq 0,$$

where $A_t = \int_0^t \frac{ds}{|X_s|^2}$, $|X| = (|X_t|, t \geq 0)$ is a Bessel process and $\mathcal{N} = (\mathcal{N}_u, u \geq 0)$ is a Poisson process independent of $|X|$. This skew-product decomposition has been extended in a nontrivial way to the multidimensional case in [3].

(b) The chaos decomposition: $(X_t)_{t\geq 0}$ decomposes as the sum of a Brownian motion $\beta$ on $\mathbb{R}$ and a purely discontinuous one dimensional martingale $\gamma$. The mixed integrals with respect to $\beta$ and $\gamma$ then span the whole $L^2$ space of $X$.

The difficulty to extend this last result to the multidimensional case is to decompose the martingale $(\eta_t)$ in the form (5) as explained above. In particular this requires a delicate study of the jumps of the Dunkl process $X$.

The organization of the paper is now presented in the form of the self-explanatory titles of the 3 next sections:

Section 2. Basic properties of the Dunkl processes.
Section 3. The martingale decomposition and a useful Itô formula.
Section 4. The chaotic representation property of the Dunkl processes.

## 2. Basic properties of the Dunkl processes.

2.1. *The intertwining operator.* The crucial tool which permits us to study the Dunkl processes is an intertwining relation between the operator $L_k$ and the Laplacian operator $\Delta$. More precisely, there exists an operator $V_k$, which turns out to be a Markovian kernel from $\mathbb{R}^d$ to itself, called the *Dunkl intertwining operator*, such that

(6) $$L_k V_k = V_k \Delta,$$

or equivalently

(7) $$P_t^{(k)} V_k = V_k P_t^{(0)},$$

where $P_t^{(k)}(x; dy) \equiv p_t^{(k)}(x, y) \, dy$ and $P_t^{(0)}(x; dy) = p_t^{(0)}(x, y) \, dy$ denote respectively the semigroups of the Dunkl process $X$ and of Brownian motion $(B_t)$.

Let us give some more details (see [9, 18]).

The intertwining operator $V_k$ is an algebraic isomorphism of the space $\mathcal{P}$ of all polynomials $P(x) = P(x^{(1)}, \ldots, x^{(d)})$ on $\mathbb{R}^d$, onto itself (the coordinates of a vector $x \in \mathbb{R}^d$ will always be written with upper indices between parentheses). Moreover, for all $i \in \mathbb{N}$, $V_k$ preserves the subspace $\mathcal{P}_i$ of all homogeneous polynomials of total degree equal to $i$. Precisely, if $\nu = (\nu^{(1)}, \ldots, \nu^{(d)}) \in \mathbb{N}^d$ is a multi-index, we consider the generalized monomial

$$(8) \qquad m_\nu(x) = V_k((\cdot)^\nu)(x), \qquad x \in \mathbb{R}^d,$$

where $(\cdot)^\nu : y \longmapsto y^\nu = (y^{(1)})^{\nu^{(1)}} \cdots (y^{(d)})^{\nu^{(d)}}$ is the usual monomial "$y$ to the power $\nu$." Then for all $i \in \mathbb{N}$, the functions $m_\nu(x)$ such that $|\nu| := \nu^{(1)} + \cdots + \nu^{(d)} = i$, constitute a basis of $\mathcal{P}_i$ and the set $\{m_\nu(x), \nu \in \mathbb{N}^d\}$ is a basis of the vector space $\mathcal{P}$. As an immediate consequence, we see that every polynomial $P(x_1, \ldots, x_l)$ with $l$ $d$-dimensional variables $x_i \in \mathbb{R}^d$ $(i = 1, \ldots, l)$ which is a linear combination of monomials $x_1^{\nu_1} \cdots x_l^{\nu_l}$ where the $\nu_i \in \mathbb{N}^d$ are multi-exponents, can also be decomposed in a linear combination of generalized monomials:

$$(9) \qquad P(x_1, \ldots, x_l) = \sum \lambda_{\nu_1, \ldots, \nu_l} m_{\nu_1}(x_1) \cdots m_{\nu_l}(x_l).$$

The intertwining operator $V_k$ has been extended by Trimèche [21] into a topological isomorphism of the space $C^\infty(\mathbb{R}^d)$ onto itself. [The space $C^\infty(\mathbb{R}^d)$ is furnished with its usual topology of the uniform convergence on compacts for the functions and all their derivatives.]

Another important tool is the *Dunkl kernel* $D_k$. It can be introduced in a concise manner as the image by $V_k$ of the usual exponential kernel. Precisely, for $x, y \in \mathbb{R}^d$, we have

$$(10) \qquad D_k(x, y) = V_k(e^{\langle \cdot, y \rangle})(x).$$

It is easily shown using (8) and the integral representation of $V_k$ given by Rösler in [19] that for all $x, y \in \mathbb{R}^d$, we have

$$(11) \qquad D_k(x, y) = \sum_{\nu \in \mathbb{N}^d} \frac{m_\nu(x)}{\nu!} y^\nu,$$

where as usual $\nu! := \nu^{(1)}! \cdots \nu^{(d)}!$.

A nontrivial result proved by Rösler [18] is that the Dunkl kernel is strictly positive on $\mathbb{R}^d \times \mathbb{R}^d$. This kernel admits a unique holomorphic extension to $\mathbb{C}^d \times \mathbb{C}^d$ and for all $z, z' \in \mathbb{C}^d$, $x \in \mathbb{R}^d$ and $\lambda \in \mathbb{C}$, we have $D_k(0, z) = 1$ and

$$(12) \quad \begin{array}{l} \text{(i) } D_k(z, z') = D_k(z', z), \\ \text{(ii) } D_k(\lambda z, z') = D_k(z, \lambda z'), \\ \text{(iii) } |\partial_z^\nu D_k(x, z)| \leq |x|^{|\nu|} \exp(|x| |Re z|), \end{array}$$

where $\partial_z^\nu$ denotes partial derivation $\nu^{(1)}$ times with respect to $z^{(1)}, \ldots, \nu^{(d)}$ times with respect to $z^{(d)}$ and $|x|$ is the Euclidean norm of $x$.

2.2. *Analytical properties of the Dunkl process.* The semigroup densities of the Dunkl process $X = (X_t)_{t \geq 0}$ with infinitesimal generator $\mathcal{L}_k$ given in (4) are expressed by the following formula obtained by Rösler [18]; see also [13]:

$$p_t^{(k)}(x,y) = \frac{1}{c_k t^{\gamma + d/2}} \exp\left(-\frac{|x|^2 + |y|^2}{2t}\right) D_k\left(\frac{x}{\sqrt{t}}, \frac{y}{\sqrt{t}}\right) \omega_k(y), \tag{13}$$

where $D_k(u,v) > 0$ is the Dunkl kernel and $\omega_k(y)$ is the *Dunkl weight function* given by

$$\omega_k(y) = \prod_{\alpha \in R_+} |\langle \alpha, y \rangle|^{2k(\alpha)}, \qquad y \in \mathbb{R}^d. \tag{14}$$

Notice that the weight function is homogeneous of degree $2\gamma$ with $\gamma$ the index of the multiplicity function of the root system

$$\gamma = \sum_{\alpha \in R_+} k(\alpha). \tag{15}$$

Moreover the constant $c_k$ is explicitly given by the integral

$$c_k = \int_{\mathbb{R}^d} e^{-|x|^2/2} \omega_k(x)\, dx. \tag{16}$$

From the explicit form (13) of $p_t^{(k)}(x,y)$, Gallardo and Yor [12] noticed the time-inversion property of the Dunkl process $X$, that is, $(tX_{1/t}, t > 0)$ is, under $\mathbb{P}_x$, a homogeneous Markov process.

Also from this explicit form (13), Rösler and Voit [20] have shown the following important result:

PROPOSITION 2.1. *The Euclidean norm $|X| = (|X_t|; t \geq 0)$ of the Dunkl process is a Bessel process of index $\alpha = \gamma + \frac{d}{2} - 1$ (i.e., of dimension $N = 2\gamma + d$), where $\gamma$ is given by* (15).

This important fact can also be deduced more directly from the explicit form of the infinitesimal generator. Indeed if in (4) we take an $SO(d)$-invariant function $u$, we can write $u(x) = F(|x|)$ for some function $F \in C_c^2(]0, +\infty[)$. Then the differences $\frac{u(x) - u(\sigma_\alpha x)}{\langle \alpha, x \rangle^2}$ vanish, $\langle \nabla u(x), \alpha \rangle = \frac{\langle \alpha, x \rangle}{|x|} F'(|x|)$ and it remains

$$\mathcal{L}_k u(x) = \frac{1}{2} F''(|x|) + \frac{d-1}{2|x|} F'(|x|) + \frac{\gamma}{|x|} F'(|x|). \tag{17}$$

We can now deduce immediately from the expression of the right-hand side of (17), that $|X|$ is a Bessel process with dimension $N = d + 2\gamma$.

We can also derive the following:

PROPOSITION 2.2. *For all $x \in \mathbb{R}^d$, under $\mathbb{P}_x$, the Dunkl process $X = (X_t, t \geq 0)$ is a martingale.*

PROOF. It is easily verified that every linear function $f_a : x \mapsto \langle a, x \rangle$ ($a \in \mathbb{R}^d$), is killed by the generator $\mathcal{L}_k$ of $X$ (i.e., $\mathcal{L}_k f_a \equiv 0$), which proves that $X$ is a local martingale. But since $|X|$ is a Bessel process, then for every $t \geq 0$, $E_x(\sup_{s \leq t} |X_s|) < +\infty$, hence $X$ is a martingale. $\square$

Let us now consider a fixed Weyl chamber $C$ of the root-system $R$ (i.e., $C$ is a connected component of $\mathbb{R}^d \setminus \bigcup_{\alpha \in R} H_\alpha$). Clearly, the space $\mathbb{R}^d/W$ of $W$-orbits in $\mathbb{R}^d$ can be identified to $\overline{C}$. Let us consider the canonical projection $\pi : \mathbb{R}^d \longrightarrow \mathbb{R}^d/W$. The process $X^W$ defined by

$$X_t^W = \pi(X_t), \tag{18}$$

is called the *W-radial part* of the Dunkl process (see [12]). It is a continuous Markov process in $\overline{C}$ with infinitesimal generator given by $\frac{1}{2} L_k^W$, where

$$L_k^W u(x) = \Delta u(x) + 2 \sum_{\alpha \in R_+} k(\alpha) \frac{\langle \nabla u(x), \alpha \rangle}{\langle \alpha, x \rangle}, \tag{19}$$

$u \in C_0^2(\overline{C}), \langle \nabla u(x), \alpha \rangle = 0$ for $x \in H_\alpha, \alpha \in R_+$. We deduce immediately from (13) that the semigroup densities of $(X_t^W)_{t>0}$ are of the form [12]:

$$p_t^W(x, y) = \frac{1}{c_k t^{\gamma + d/2}} \exp\left(-\frac{|x|^2 + |y|^2}{2t}\right) D_k^W\left(\frac{x}{\sqrt{t}}, \frac{y}{\sqrt{t}}\right) \omega_k(y), \qquad x, y \in C,$$

where

$$D_k^W(u, v) = \sum_{w \in W} D_k(u, wv).$$

It is interesting to notice that the Brownian motion process in a Weyl chamber studied by Biane, Bougerol and O'Connell in [2] is a particular $W$-radial Dunkl process corresponding to the multiplicity function $k(\alpha) \equiv 1$ (see [12]).

### 3. The martingale decomposition and a useful Itô formula.

3.1. *The jumps of the Dunkl process.* Let us first recall that the Lévy kernel $N(x, dy)$ of a homogeneous Markov process with a semigroup $(P_t)_{t \geq 0}$ of transition operators, is determined, for all $x \in \mathbb{R}^d$, by

$$\lim_{t \to 0} \frac{P_t f(x)}{t} = \int_{\mathbb{R}^d} N(x, dy) f(y),$$

for $f$ a function in the domain of the infinitesimal generator which vanishes in a neighborhood of $x$ [16]. For the Dunkl process $X$, and with the notation introduced in Section 1, we have the following result:

PROPOSITION 3.1. *The Lévy kernel of $X$ has the following form*

$$
(20) \quad N(x,dy) = \begin{cases} \displaystyle\sum_{\alpha \in R_+} k(\alpha)\frac{\delta_{\sigma_\alpha x}(dy)}{\langle \alpha, x\rangle^2}, & \text{if } x \notin \bigcup_{\alpha \in R_+} H_\alpha, \\ \displaystyle\sum_{\substack{\alpha \notin S \\ \alpha \in R_+}} k(\alpha)\frac{\delta_{\sigma_\alpha x}(dy)}{\langle \alpha, x\rangle^2}, & \text{if } x \in \bigcap_{\alpha \in S} H_\alpha, \\ 0, & \text{if } x = 0, \end{cases}
$$

*where $S$ denotes any subset of the root system $R_+$ and $\delta_z(dy)$ denotes the Dirac measure at point $z \in \mathbb{R}^d$.*

PROOF. For $f \in C^2(\mathbb{R}^d) \cap \text{Dom}(\mathcal{L}_k)$ and $f = 0$ in a neighborhood of $x$, it suffices to write $\frac{1}{t} P_t f(x) = \mathcal{L}_k f(x) + \varepsilon_t$, where $\varepsilon_t \to 0$ if $t \to 0$ and to use (4). □

Now consider the function $f$ defined on $\mathbb{R}^d \times \mathbb{R}^d$ as $f(x,y) = \mathbf{1}_{(\text{supp}N(x,dz))^c}(y)$, that is, taking the value 1 if $y$ is not in the set $\{\sigma_\alpha x, \alpha \in R_+\}$ and the value 0 otherwise. We have the following result, which is, clearly, a very particular illustration of the general fact that the Lévy kernel of a Markov process $X$ "rules" its jumps.

PROPOSITION 3.2. *For all $t > 0$ and all $x \in \mathbb{R}^d$, we have*

$$
(21) \quad \sum_{s \leq t} f(X_{s-}, X_s)\mathbf{1}_{(X_s \neq X_{s-})} = 0, \qquad \mathbb{P}_x\text{-a.s.}
$$

PROOF. The positive discontinuous additive functional on the left-hand side of (21) can be compensated by the process

$$
(22) \quad \int_0^t ds \int_{\mathbb{R}^d} N(X_{s-}, dy) f(X_{s-}, y),
$$

which is identically zero according to the definition of $f$ (see [16], page 154 or [22]). The result then clearly follows. □

REMARK. The result (21) says that when there is a jump of the process at (a random) time $s$, that is, $X_s \neq X_{s-}$, then, almost surely, there exists (a random) $\alpha \in R_+$ such that $X_s = \sigma_\alpha X_{s-}$. We will say that the jump occurs in the direction $\alpha$. In that case, according to (3),

$$
\Delta X_s := X_s - X_{s-} = -\langle \alpha, X_{s-}\rangle \alpha.
$$

Clearly from the form of the Lévy kernel and the above discussion, if $k(\alpha) = 0$ for some $\alpha \in R_+$, no jump can occur in the direction $\alpha$.

We will now show that the sum over any time interval of the amplitudes of the jumps is finite. Precisely:

PROPOSITION 3.3. *For every $t > 0$ and $x \in \mathbb{R}^d \setminus \{0\}$,*

$$\sum_{s \leq t} |\Delta X_s| < +\infty, \qquad \mathbb{P}_x\text{-}a.s. \tag{23}$$

PROOF. According to the above remark, we can write

$$\sum_{s \leq t} |\Delta X_s| = \sum_{\alpha \in R_+} \sum_{s \leq t} f_\alpha(X_{s-}, X_s),$$

where

$$f_\alpha(x, y) = \sqrt{2} |\langle \alpha, x \rangle| \mathbf{1}_{(y = \sigma_\alpha x \neq x)}.$$

But the positive discontinuous functional $\sum_{s \leq t} f_\alpha(X_{s-}, X_s)$ can be compensated (with respect to quantities which appear via compensation, we keep writing $X_{s-}$, although $X_s$ would also be correct since $\{s : X_s \neq X_{s-}\}$ is Lebesgue negligible) by the process

$$\int_0^t ds \int_{\mathbb{R}^d} N(X_{s-}, dy) f_\alpha(X_{s-}, y) = \int_0^t ds \int_{\mathbb{R}^d} \sum_{\beta \in R_+} k(\beta) \frac{\delta_{\sigma_\beta X_{s-}}(dy)}{\langle \beta, X_{s-} \rangle^2} f_\alpha(X_{s-}, y)$$

$$= \sqrt{2} k(\alpha) \int_0^t \frac{ds}{|\langle \alpha, X_{s-} \rangle|}.$$

Thus, it will suffice to prove the following lemma:

LEMMA 3.4. *For every $\alpha \in R_+$ such that $k(\alpha) > 0$, for all $x \in \mathbb{R}^d \setminus \{0\}$ and $t > 0$, we have*

$$E_x \left( \int_0^t \frac{ds}{|\langle \alpha, X_{s-} \rangle|} \right) < +\infty. \tag{24}$$

Indeed, this implies that the sum of jumps (23) has a finite expectation and the proposition follows. □

REMARK. This result may be compared with the finiteness of

$$E_x \left( \int_0^t \frac{ds}{|X_{s-}|} \right) < +\infty, \tag{25}$$

which occurs as $|X_s|$ is a Bessel process with dimension $N = 2\gamma + d > 1$. However (24) is finer than (25) because $|\langle \alpha, X_{s-} \rangle| \leq \sqrt{2} |X_{s-}|$.

PROOF OF LEMMA 3.4. By the expression (13) of the transition densities, we can write

$$
\begin{aligned}
E_x&\left(\int_0^t \frac{ds}{|\langle \alpha, X_{s-}\rangle|}\right)\\
&= \int_0^t ds\, \frac{e^{-|x|^2/(2s)}}{c_k s^{\gamma+d/2}} \int_{\mathbb{R}^d} e^{-|y|^2/(2s)} D_k\left(\frac{x}{\sqrt{s}}, \frac{y}{\sqrt{s}}\right) \frac{\omega_k(y)}{|\langle \alpha, y\rangle|}\, dy.
\end{aligned}
\tag{26}
$$

But by the homogeneity of the weight function $\omega_k$ and the fact that $D_k(\frac{x}{\sqrt{s}}, \frac{y}{\sqrt{s}}) = D_k(x, \frac{y}{s})$ [see (12)(ii)], we can perform the change of variables $y = sz$ in the right-hand side of (26) and we see that the quantity (26) is equal to

$$
\frac{1}{c_k} \int_0^t e^{-|x|^2/(2s)} s^{\gamma+(d/2)-1}\, ds \int_{\mathbb{R}^d} e^{-s|z|^2/2} D_k(x, z) \frac{\omega_k(z)}{|\langle \alpha, z\rangle|}\, dz.
\tag{27}
$$

Now using $D_k(x, z) \le e^{|x||z|}$ [see (12)(iii)], the explicit expression of $\omega_k$ (14) and the Cauchy–Schwarz inequality, we can see that the integral on $\mathbb{R}^d$ with respect to $z$ in the right-hand side of (27), is bounded by

$$
C \int_{\mathbb{R}^d} e^{-s|z|^2/2} e^{|x||z|} |z|^{2\gamma - 2k(\alpha)} |\langle \alpha, z\rangle|^{2k(\alpha)-1}\, dz,
\tag{28}
$$

with $C$ a positive constant. In order to get an estimate of the integral (28), we can suppose, after performing, if necessary, a suitable rotation of the coordinate axes, that $\alpha = \sqrt{2} e_d$, where $e_d$ is the $d$th vector of the canonical basis of $\mathbb{R}^d$. Then if we use generalized polar coordinates in $\mathbb{R}^d$ involving angles $\theta_1 \in [0, \pi/2]$ and $\theta_i \in [-\pi/2, \pi/2]$ for $i = 1, \ldots, d-1$ we have $\langle \alpha, z\rangle = \sqrt{2}\rho \sin\theta_{d-1}$ and the integral (28) is equal to

$$
C \int_0^\infty e^{-s\rho^2/2} e^{|x|\rho} \rho^{2\gamma+d-2}\, d\rho \int_{S^{d-1}} |\sin\theta_{d-1}|^{2k(\alpha)-1}\, d\sigma,
\tag{29}
$$

where $S^{d-1}$ is the unit sphere and

$$
d\sigma = (\cos\theta_{d-1})^{d-2}(\cos\theta_{d-2})^{d-3} \cdots (\cos\theta_2)\, d\theta_1\, d\theta_2 \cdots d\theta_{d-1}
$$

is the uniform measure on $S^{d-1}$. The integral over $S^{d-1}$ in (29) converges because $k(\alpha) > 0$. Then using (29), (28) and interchanging the summations, we see that the integral (27) is smaller than

$$
C \int_0^\infty H(\rho) e^{|x|\rho} \rho^{2\gamma+d-2}\, d\rho,
\tag{30}
$$

where $C$ is another positive constant and

$$
H(\rho) = \int_0^t s^{\gamma+(d/2)-1} e^{-s\rho^2/2} e^{-|x|^2/(2s)}\, ds.
\tag{31}
$$

It remains to prove that the integral (30) is convergent. In order to find an equivalent of $H(\rho)$ as $\rho \to +\infty$, set $y = \rho^2/2$, then $H(\sqrt{2y}) := \varpi(y)$ is the Laplace transform of the distribution function $U$ of a measure concentrated on $[0, t]$ and given for $s \in [0, +\infty[$ by

$$U(s) = \int_0^{s \wedge t} v^{\gamma + (d/2) - 1} e^{-|x|^2/(2v)} \, dv. \tag{32}$$

Clearly, we can apply the Tauberian theorem (see [10], page 443) and we get $\varpi(y) \sim U(1/y)$ if $y \to +\infty$. Then from (32) with $s = 1/y$ and changing $v$ in $(u+y)^{-1}$, we obtain

$$\begin{aligned}
U(1/y) &= \int_0^\infty (y+u)^{-(\gamma + (d/2) + 1)} e^{-(y+u)|x|^2/2} \, du \\
&\sim 2|x|^{-2} y^{-(\gamma + (d/2) + 1)} e^{-y|x|^2/2},
\end{aligned} \tag{33}$$

as $y \to +\infty$. Using (33) and returning to $H$, we obtain the equivalent as $\rho \to +\infty$:

$$H(\rho) \sim 2|x|^{-2} (\rho^2/2)^{-(\gamma + (d/2) + 1)} e^{-\rho^2 |x|^2/4}, \tag{34}$$

we immediately see that the integral (30) is convergent and the lemma is proved. $\square$

With the same ideas using the Lévy kernel, we can give some information on the number

$$N_t = \sum_{s \leq t} \sum_{\alpha \in R_+} \mathbf{1}_{(X_s = \sigma_\alpha X_{s-} \neq X_{s-})} \tag{35}$$

of jumps of the process $X$ on the time interval $[0, t]$; more precisely, we have the following result:

PROPOSITION 3.5. *If the multiplicity function of the root system satisfies $k(\alpha) > 1/2$ for all $\alpha \in R_+$, then for all $x \in \mathbb{R}^d \setminus \{0\}$ and $t > 0$, we have*

$$N_t < +\infty, \qquad \mathbb{P}_x\text{-}a.s.$$

PROOF. The positive discontinuous functional given in (35) may be easily compensated by the process

$$\sum_{\alpha \in R_+} k(\alpha) \int_0^t \frac{ds}{\langle \alpha, X_{s-} \rangle^2}$$

which then has the same expectation as $N_t$. This expectation is finite if for all $\alpha \in R_+$, we have

$$E_x \left( \int_0^t \frac{ds}{\langle \alpha, X_{s-} \rangle^2} \right) < +\infty. \tag{36}$$

This is obtained by the same method used in the proof of Lemma 3.4, replacing $|\langle \alpha, y \rangle|$ by $\langle \alpha, y \rangle^2$ in (26), (27) and (28). In (29) the exponent of $|\sin \theta_{d-1}|$ has to be replaced by $2k(\alpha) - 2$; then the integral over the sphere $S^{d-1}$ is convergent if $k(\alpha) > 1/2$. The end of the proof is unchanged. □

REMARK. The expectation (36) multiplied by the factor $k(\alpha)$, represents the mean number of the jumps of the process $X$ in the direction of the root $\alpha$.

3.2. *The martingale decomposition.* We now state the main result of this section.

THEOREM 1. *For all $x \in \mathbb{R}^d \setminus \{0\}$, under $P_x$, we can write the martingale decomposition of $(X_t)$ into its continuous and purely discontinuous parts, in the following form*

$$(37) \qquad X_t = x + B_t + \sum_{\alpha \in R_+} \sqrt{k(\alpha)} M_t^\alpha \alpha,$$

*where $(B_t, t \geq 0)$ is a d-dimensional Brownian motion and $(M_t^\alpha, t \geq 0)_{\alpha \in R_+}$ is a family of "normal" martingales, in the sense of Meyer [17]:*

$$(38) \qquad \langle M^\alpha, M^\alpha \rangle_t = t$$

*and*

$$(39) \qquad [M^\alpha, M^\beta]_t = 0 \qquad \text{for } \alpha \neq \beta.$$

*Moreover, $(M_t^\alpha, t \geq 0)$ has paths of finite variation, and may be written as the compensated sum of its jumps:*

$$(40) \qquad M_t^\alpha = \sum_{s \leq t} \frac{-\langle \alpha, X_{s-} \rangle}{\sqrt{k(\alpha)}} \mathbf{1}_{\{X_s = \sigma_\alpha(X_{s-}) \neq X_{s-}\}} + \int_0^t \frac{\sqrt{k(\alpha)}\, ds}{\langle \alpha, X_{s-} \rangle}.$$

PROOF. (i) (the continuous part): For every function $\varphi \in C^2(\mathbb{R}^d)$, the general Itô's formula writes (we are using Einstein's convention, i.e., simple or double summation is understood with respect to indices both in upper and lower positions)

$$(41) \quad \begin{aligned} \varphi(X_t) &= \varphi(X_0) + \int_0^t \partial_i \varphi(X_{s-})\, dX_s^{(i)} + \tfrac{1}{2} \int_0^t \partial_{ij}^2 \varphi(X_s)\, d[X^{(i)}, X^{(j)}]_s^c \\ &\quad + \sum_{s \leq t}(\varphi(X_s) - \varphi(X_{s-}) - \partial_i \varphi(X_{s-}) \Delta X_s^{(i)}). \end{aligned}$$

But the processes

$$(42) \qquad \varphi(X_t) - \int_0^t \tfrac{1}{2} L_k \varphi(X_s)\, ds$$

and

$$\varphi(X_0) + \int_0^t \partial_i \varphi(X_{s-}) \, dX_s^{(i)} \tag{43}$$

are local martingales as well as their difference which by (41) is equal to

$$\sum_{s \leq t} (\varphi(X_s) - \varphi(X_{s-}) - \partial_i \varphi(X_{s-}) \Delta X_s^{(i)})$$

$$- \int_0^t \tfrac{1}{2} L_k \varphi(X_s) \, ds + \tfrac{1}{2} \int_0^t \partial_{ij}^2 \varphi(X_s) \, d[X^{(i)}, X^{(j)}]_s^c.$$

This immediately shows that the process

$$\int_0^t \tfrac{1}{2} L_k \varphi(X_s) \, ds - \tfrac{1}{2} \int_0^t \partial_{ij}^2 \varphi(X_s) \, d[X^{(i)}, X^{(j)}]_s^c \tag{44}$$

is the compensator of the discontinuous functional

$$\sum_{s \leq t} (\varphi(X_s) - \varphi(X_{s-}) - \partial_i \varphi(X_{s-}) \Delta X_s^{(i)}). \tag{45}$$

But by the theory of the Lévy kernel (see, e.g., [16]), this compensator is also equal to

$$\int_0^t ds \int_{\mathbb{R}^d} N(X_{s-}, dy)(\varphi(y) - \varphi(X_{s-}) - \partial_i \varphi(X_{s-})(y - X_{s-})^{(i)}). \tag{46}$$

Using the expression (20) we easily see that (46) is equal to

$$\int_0^t ds \left\{ \sum_{\alpha \in R_+} k(\alpha) \left[ \frac{\langle \nabla \varphi(X_{s-}), \alpha \rangle}{\langle \alpha, X_{s-} \rangle} + \frac{\varphi(\sigma_\alpha X_{s-}) - \varphi(X_{s-})}{\langle \alpha, X_{s-} \rangle^2} \right] \right\}. \tag{47}$$

Finally from the equality between (44) and (47) and the expression (4) of $L_k$, we deduce that for every function $\varphi \in C^2(\mathbb{R}^d)$, we have

$$\int_0^t \Delta \varphi(X_s) \, ds = \int_0^t \partial_{ij}^2 \varphi(X_s) \, d[X^{(i)}, X^{(j)}]_s^c. \tag{48}$$

Applying (48) to the functions $\varphi(x) = x^{(i)} x^{(j)}$, $i,j = 1, \ldots, d$, we get $[X^{(i)}, X^{(j)}]_t^c = 0$ if $i \neq j$ and $[X^{(i)}, X^{(i)}]_t^c = t$ for all $i$. But as $[X^{(i)}, X^{(j)}]_t^c = \langle (X^{(i)})^c, (X^{(j)})^c \rangle_t$, Lévy's theorem implies that $X^c$ is a $d$-dimensional Brownian motion.

(ii) (the purely discontinuous part): By Proposition 3.3, the process $\sum_{s \leq t} \Delta X_s$ is well defined and by the remark following Proposition 3.2, we can write

$$\sum_{s \leq t} \Delta X_s = \sum_{\alpha \in \mathbb{R}_+} S_t^\alpha \alpha \tag{49}$$

where for all $\alpha \in R_+$,

$$
S_t^\alpha = \sum_{s \leq t} -\langle \alpha, X_{s-}\rangle \mathbf{1}_{(X_s = \sigma_\alpha X_{s-} \neq X_{s-})}, \tag{50}
$$

is the sum of the jumps of $X$ in the direction of the root $\alpha$. With the same method used in the proof of Proposition 3.3, we show that the functional (50), is compensated by the process $C^\alpha = \{C_t^\alpha = -k(\alpha) \int_0^t \frac{ds}{\langle \alpha, X_{s-}\rangle}, t \geq 0\}$, which is well defined and has finite variation by Lemma 3.4. This implies that the continuous and finite variation process $C = \sum_{\alpha \in R_+} C^\alpha$ is the compensator of the sum of jumps (49). Now the decomposition

$$
X_t = X_t - \sum_{s \leq t} \Delta X_s + C_t + \sum_{s \leq t} \Delta X_s - C_t \tag{51}
$$

shows that the continuous martingale part of $(X_t)$ (which is unique) is equal to $X_t - \sum_{s \leq t} \Delta X_s + C_t$ and $X_t^d = \sum_{s \leq t} \Delta X_s - C_t$ which is a finite variation martingale, hence purely discontinuous, is its purely discontinuous martingale part. Clearly, if we let $M_t^\alpha = (k(\alpha))^{-1/2}(S_t^\alpha - C_t^\alpha)$, we obtain (40) and (51) can be rewritten in the form (37). Finally, as $M^\alpha$ is a finite variation martingale, we have

$$
\begin{aligned}
[M^\alpha, M^\alpha]_t &= \sum_{s \leq t} (\Delta M_s^\alpha)^2 \\
&= (k(\alpha))^{-1} \sum_{s \leq t} \langle \alpha, X_{s-}\rangle^2 \mathbf{1}_{(X_s = \sigma_\alpha X_{s-} \neq X_{s-})},
\end{aligned} \tag{52}
$$

and this functional (52) can be compensated by a process of the form given in (22) which, in this case, is equal to $t$. This proves assertion (38). Assertion (39) follows immediately from the fact that two jumps of the process respectively in the directions $\alpha$ and $\beta$ cannot occur simultaneously. The theorem is then entirely proved. $\square$

COROLLARY 3.6 (An Itô formula). *Under the notation of Theorem 1, for any function $f \in C^{2,1}(\mathbb{R}^d \times \mathbb{R})$, the semimartingale $f(X_t, t)$ can be decomposed in the following form:*

$$
\begin{aligned}
f(X_t, t) = f(X_0, 0) &+ \int_0^t \langle \nabla_x f(X_s, s), dB_s\rangle \\
&+ \sum_{\alpha \in R_+} \int_0^t \frac{f(X_{s-}, s) - f(\sigma_\alpha X_{s-}, s)}{\langle \alpha, X_{s-}\rangle} \sqrt{k(\alpha)}\, dM_s^\alpha \\
&+ \int_0^t ds \left(\frac{\partial f}{\partial s} + \frac{1}{2} L_k f\right)(X_s, s).
\end{aligned} \tag{53}
$$

PROOF. The usual Itô formula writes

$$f(X_t, t) = f(X_0, 0) + \int_0^t \partial_i f(X_{s-}, s) \, dX_s^{(i)}$$

(54)
$$+ \int_0^t \frac{\partial f}{\partial s}(X_s, s) \, ds + \frac{1}{2} \int_0^t \partial_{ij}^2 f(X_s, s) \, d[X^{(i)}, X^{(j)}]_s^c$$

$$+ \sum_{s \leq t}(f(X_s, s) - f(X_{s-}, s) - \partial_i f(X_{s-}, s)\Delta X_s^{(i)}).$$

But the process

(55) $$M_t^f \stackrel{\text{def}}{=} f(X_t, t) - \int_0^t \left(\frac{\partial f}{\partial s} + \frac{1}{2}L_k f\right)(X_s, s) \, ds$$

is a local martingale. Then, taking into account the fact that $d[X^{(i)}, X^{(j)}]_s^c = \delta_{ij} \, ds$, and that $f(X_0, 0) + \int_0^t \partial_i f(X_{s-}, s) \, dX_s^{(i)}$ is a local martingale, if we substract the term $\int_0^t (\frac{\partial f}{\partial s} + \frac{1}{2}L_k f)(X_s, s) \, ds$ from the left-hand side of (54), we obtain that

(56)
$$\sum_{s \leq t}(f(X_s, s) - f(X_{s-}, s) - \partial_i f(X_{s-}, s)\Delta X_s^{(i)})$$
$$- \sum_{\alpha \in R_+} k(\alpha) \int_0^t \left[\frac{\langle \nabla_x f(X_s, s), \alpha \rangle}{\langle \alpha, X_{s-} \rangle} - \frac{f(X_{s-}, s) - f(\sigma_\alpha X_{s-}, s)}{\langle \alpha, X_{s-} \rangle^2}\right] ds$$

is a local martingale. From the structure of the jumps given by Theorem 1, the first sum in (56) is equal to

(57) $$\sum_{s \leq t} \sum_{\alpha \in R_+} \left(\frac{f(X_{s-}, s) - f(\sigma_\alpha X_{s-}, s)}{\langle \alpha, X_{s-} \rangle} - \langle \nabla_x f(X_{s-}, s), \alpha \rangle\right)\sqrt{k(\alpha)}\Delta M_s^\alpha$$

and we can then easily show that the local martingale (56) is equal to

(58) $$\int_0^t \sum_{\alpha \in R_+} \left(\frac{f(X_{s-}, s) - f(\sigma_\alpha X_{s-}, s)}{\langle \alpha, X_{s-} \rangle} - \langle \nabla_x f(X_{s-}, s), \alpha \rangle\right)\sqrt{k(\alpha)} \, dM_s^\alpha.$$

If in formula (54) we replace the expression of (57) by the value obtained from (58) and we rearrange the sum, then we immediately get (53). □

Note that Corollary 3.6 provides the explicit decomposition of the local martingale (55) into its continuous and purely discontinuous parts. In particular, formula (53) shows that our study (but not our method!) falls precisely in the general scheme adopted by Dellacherie–Maisonneuve–Meyer ([4], page 259 and sequel; see, e.g., formula (62.5) there) to prove the chaotic representation property for a general class of Markov processes, following Biane [1] and Isobe and Sato [14]. Moreover, Corollary 3.6 leads to the following predictable representation property:

COROLLARY 3.7. *For any $x \in \mathbb{R}^d \setminus \{0\}$, under $\mathbb{P}_x$, every local martingale in the filtration of $X$, may be written as*

$$C + \int_0^t \langle \mu_c(s), dB_s \rangle + \sum_{\alpha \in R_+} \int_0^t \mu_\alpha(s) \, dM_s^\alpha,$$

*where $C$ is a constant, and $\mu_c$ and $\mu_\alpha$ ($\alpha \in R_+$) are predictable processes such that*

$$\int_0^t |\mu_c(s)|^2 \, ds < +\infty, \qquad \int_0^t (\mu_\alpha(s))^2 \, ds < +\infty, \qquad \mathbb{P}_x\text{-}a.s.$$

PROOF. From [5] (Théorème 25, page 246) it is known that the martingales $M_t^f$ given by (55), for $f \in C^{2,1}(\mathbb{R}^d \times \mathbb{R})$, generate in the sense of Kunita–Watanabe, the space of locally square integrable martingales. However from Corollary 3.6, $M_t^f$ is itself the sum of stochastic integrals with respect to $dB_s$ and $dM_s^\alpha$. The result follows. □

### 4. The chaos decomposition.

4.1. *The Dunkl–Hermite space–time polynomials.* Let $B = (B_t)_{t \geq 0}$ be a $d$-dimensional Brownian motion and let $\nu \in \mathbb{N}^d$ be a multi-index. Then for $0 \leq s \leq t$, we have

(59) $$E(B_t^\nu | \mathcal{F}_s) = \mathcal{H}_\nu(B_s, s-t),$$

where $\mathcal{H}_\nu(x,t)$ ($x \in \mathbb{R}^d, t \in \mathbb{R}$) is the $\nu$th space–time polynomial given by

(60) $$\mathcal{H}_\nu(x,t) = \prod_1^d \mathcal{H}_{\nu^{(i)}}(x^{(i)}, t),$$

with $\nu = (\nu^{(1)}, \ldots, \nu^{(d)})$, $x = (x^{(1)}, \ldots, x^{(d)})$ and for $n \in \mathbb{N}$, $\mathcal{H}_n(x,t)$ is the classical $n$th space–time Hermite polynomial (see [23]). A key property of $\mathcal{H}_\nu$ is its space–time harmonicity, that is,

(61) $$\left(\frac{\partial}{\partial t} + \frac{1}{2}\Delta\right)\mathcal{H}_\nu(x,t) = 0.$$

We will now establish the analog of (59) and (61) for the Dunkl process $X$. We use the notation of Section 2.1 [in particular (8)]:

PROPOSITION 4.1. *For all $\nu \in \mathbb{N}^d$ and $s \leq t$, we have*

(62) $$E(m_\nu(X_t)|\mathcal{F}_s) = Q_\nu(X_s, s-t),$$

*where $Q_\nu(x,t)$ is a space–time polynomial on $\mathbb{R}^d \times \mathbb{R}$ of total degree $|\nu|$ in the space variable. Moreover this polynomial is space–time harmonic, that is*

(63) $$\left(\frac{\partial}{\partial t} + \frac{1}{2}L_k\right)Q_\nu(x,t) = 0$$

*(where $L_k$ is acting on the space variable).*

PROOF. The Markov property yields $E(m_\nu(X_t)|\mathcal{F}_s) = P^{(k)}_{t-s}(m_\nu)(X_s)$ a.s. But for all $y \in \mathbb{R}^d$, by (7) and (8), we can write [for an operator $T$ and a function $f$ we write abusively $T(f(y))$ for $T(f)(y)$]

$$P_t^{(k)}(m_\nu(y)) = P_t^{(k)} V_k(y^\nu)$$
$$= V_k P_t^{(0)}(y^\nu)$$
$$= V_k \mathcal{H}_\nu(y, -t).$$

As the operator $V_k$ (which acts on the space variable) preserves polynomials and degrees, $V_k \mathcal{H}_\nu(x, t) := Q_\nu(x, t)$ is a polynomial of degree $|\nu|$. This proves (62). Moreover using (6) and the commutation of the operators $V_k$ and $\frac{\partial}{\partial t}$, we have

$$\left(\frac{\partial}{\partial t} + \frac{1}{2} L_k\right) Q_\nu(y, t) = \frac{\partial}{\partial t} V_k \mathcal{H}_\nu(y, t) + \frac{1}{2} V_k \Delta \mathcal{H}_\nu(y, t)$$
$$= V_k \left(\frac{\partial}{\partial t} \mathcal{H}_\nu(y, t) + \frac{1}{2} \Delta \mathcal{H}_\nu(y, t)\right)$$
$$= 0.$$

Then (63) follows. □

COROLLARY 4.2. *Under the notation of Proposition* 4.1 *and Theorem* 1, *for all $s' \leq s \leq t$, we have*

(64)
$$Q_\nu(X_s, s - t) = Q_\nu(X_{s'}, s' - t) + \sum_{i=1}^{d} \int_{s'}^{s} Q_{\nu,c}^{(i)}(X_{u-}, u - t) \, dB_u^{(i)}$$
$$+ \sum_{\alpha \in R_+} \int_{s'}^{s} Q_{\nu,\delta}^\alpha(X_{u-}, u - t) \, dM_u^\alpha,$$

*where $Q_{\nu,c}^{(i)}(x, t)$ and $Q_{\nu,\delta}^\alpha(x, t)$ are space–time polynomials of degree $\leq |\nu| - 1$ in the space variable which are given by*

(65)
$$Q_{\nu,c}^{(i)}(x, t) = \frac{\partial}{\partial x^{(i)}} Q_\nu(x, t),$$
$$Q_{\nu,\delta}^\alpha(x, t) = \sqrt{k(\alpha)} \frac{Q_\nu(x, t) - Q_\nu(\sigma_\alpha x, t)}{\langle \alpha, x \rangle}.$$

PROOF. If we apply the Itô formula (53) to the space–time harmonic function $f(x, u) = Q_\nu(x, u - t)$ (with $t$ fixed) and we substitute $f(X_s, s) - f(X_{s'}, s')$, we obtain immediately (65). Moreover, as $\langle \alpha, x \rangle = 0$ implies $Q_\nu(x, t) - Q_\nu(\sigma_\alpha x, t) = 0$, $Q_{\nu,\delta}^\alpha(x, t)$ is clearly a polynomial of degree $\leq |\nu| - 1$ in the space variable and the result follows. □

4.2. *The Wiener chaos decomposition.* Let us consider the integer $N = d + \operatorname{card}(R_+)$. In order to describe the $n$th Wiener chaos, we define a numbering of the roots $\alpha \in R_+$ of the form $\alpha_{d+1}, \ldots, \alpha_N$ and for any integer $\varepsilon \in [1, N]$, we denote by $Z^\varepsilon = (Z^\varepsilon_t)_{t \geq 0}$, the real-valued martingale given by

$$(66) \qquad Z^\varepsilon_t = \begin{cases} B^{(\varepsilon)}_t, & \text{if } \varepsilon \in [1, d], \\ M^{\alpha_\varepsilon}_t, & \text{if } \varepsilon \in [d+1, N]. \end{cases}$$

The Wiener chaos $C_n(X)$ of order $n \geq 1$ is then defined as the subspace of $L^2(\mathcal{F}_\infty)$ generated by the "mixed" stochastic integrals

$$(67) \qquad \int_0^\infty dZ^{\varepsilon_1}_{u_1} \int_0^{(u_1)-} dZ^{\varepsilon_2}_{u_2} \cdots \int_0^{(u_{n-1})-} dZ^{\varepsilon_n}_{u_n} f_n(u_1, \ldots, u_n),$$

for $(\varepsilon_1, \ldots, \varepsilon_n) \in [1, N]^n$ and $f_n : \Delta_n \to \mathbb{R}$ a measurable function such that:

$$\int_{\Delta_n} du_1 \cdots du_n \, f_n^2(u_1, \ldots, u_n) < \infty,$$

where

$$\Delta_n = \{(u_1, \ldots, u_n); 0 < u_n < u_{n-1} < \cdots < u_1\}.$$

If $n = 0$, let $C_0(X)$ be the subspace of the constants. With these notation, the main result of our paper is the following:

THEOREM 2. $L^2(\mathcal{F}_\infty) = \bigoplus_{n=0}^\infty C_n(X)$.

PROOF. By Proposition 2.1, the coordinates of $(X_t)_{t \geq 0}$ have all exponential moments. It is then well known that the subspace generated by the random variables of the form $P(X_{t_1}, \ldots, X_{t_l})$, for $l \in \mathbb{N} \setminus \{0\}$, $0 \leq t_1 < t_2 < \cdots < t_l$ and $P(x_1, \ldots, x_l)$ a polynomial on $(\mathbb{R}^d)^l$, is dense in $L^2(\mathcal{F}_\infty)$. Consequently, by (9), the random variables

$$(68) \qquad \Lambda = \prod_{j=1}^l m_{\nu_j}(X_{t_j}),$$

with $l \in \mathbb{N} \setminus \{0\}$, $0 \leq t_1 < t_2 < \cdots < t_l$, $\nu_1, \ldots, \nu_l \in \mathbb{N}^d$, form a total subset of $L^2(\mathcal{F}_\infty)$. The theorem will then be proved if we can show that every random variable $\Lambda$ of the form (68) can be written as a finite sum of stochastic integrals of the type (67). But by (62), we have $m_{\nu_l}(X_{t_l}) = Q_{\nu_l}(X_{t_l}, 0)$. Then if we apply (64) with $s = t = t_l$ and $s' = t_{l-1}$, we obtain

$$(69) \qquad \begin{aligned} m_{\nu_l}(X_{t_l}) &= Q_{\nu_l}(X_{t_{l-1}}, t_{l-1} - t_l) \\ &\quad + \sum_{i=1}^d \int_{t_{l-1}}^{t_l} Q^{(i)}_{\nu_l, c}(X_{u-}, u - t_l) \, dB^{(i)}_u \\ &\quad + \sum_{\alpha \in R_+} \int_{t_{l-1}}^{t_l} Q^\alpha_{\nu_l, \delta}(X_{u-}, u - t_l) \, dM^\alpha_u. \end{aligned}$$

The random variable $\Lambda$ can then be written as a sum of the following terms, where the times $t_1, t_2, \ldots, t_{l-1}, u, t_l$ are ranked in increasing order:

(i) A random variable

$$\Lambda_1 = Q_{\nu_l}(X_{t_{l-1}}, t_{l-1} - t_l) \prod_{j=1}^{l-1} m_{\nu_j}(X_{t_j}),$$

which depends on $l-1$ times;

(ii) A sum of stochastic integrals with respect to $dB_u^{(i)}$ ($i=1,\ldots,d$) with integrands

$$\Lambda_2^{(i)}(u) = Q_{\nu_l,c}^{(i)}(X_{u-}, u - t_l) \prod_{j=1}^{l-1} m_{\nu_j}(X_{t_j}),$$

which is a polynomial depending on the times $t_1, \ldots, t_{l-1}, u$ but has degree $\leq |\nu_1| + \cdots + |\nu_l| - 1$ (in the space variable); and finally:

(iii) a sum of card $R_+ = N - d$ stochastic integrals with respect to $dM_u^\alpha$ ($\alpha \in R_+$) with integrands,

$$\Lambda_3^\alpha(u) = Q_{\nu_l,\delta}^\alpha(X_{u-}, u - t_l) \prod_{j=1}^{l-1} m_{\nu_j}(X_{t_j}),$$

which depends on the times $t_1, \ldots, t_{l-1}, u$ but where the polynomial $Q_{\nu_l,\delta}^\alpha$ is of degree $\leq |\nu_l| - 1$.

The chaotic decomposition then follows by a decreasing induction on $l$ and the total degree $|\nu_1| + \cdots + |\nu_l|$ with the same method used in [23], Chapter 15, for the chaotic decomposition of the Azéma–Emery martingales.

**Acknowledgments.** We thank Michel Emery for several fruitful discussions related to the martingale decompositions of the Dunkl processes and for pointing out some similarities with multidimensional Azéma martingales. We also thank the referees for their useful remarks.

Laboratoire de Mathématiques
  et Physique Théorique
CNRS-UMR 6083
Université de Tours
Campus de Grandmont
37200 Tours
France
E-mail: gallardo@univ-tours.fr

Laboratoire de Probabilités
  et Modèles Aléatoires
CNRS-UMR 7599
Université de Paris 6 et IUF
4 Place Jussieu
75252 Paris
France
E-mail: deaproba@proba.jussieu.fr